\documentclass[a4paper,11pt]{article}
\usepackage[verbose,left=29mm,right=29mm,top=33mm,bottom=25mm]{geometry}
\usepackage[english]{babel}
\usepackage[utf8]{inputenc}
\usepackage{graphicx}
\usepackage{amsmath,amsfonts,amscd,bezier,amsthm,amssymb}
\usepackage[T1]{fontenc}
\usepackage{todonotes}
\usepackage{natbib}
\usepackage{longtable}
\usepackage{amsmath}
\usepackage{amsthm}
\usepackage{eucal}
\usepackage{amssymb}
\usepackage{multirow}
\usepackage{mathrsfs}
\usepackage{algorithm}
\usepackage{algorithmic}
\usepackage{booktabs} % To thicken table lines
\usepackage{rotating}
\usepackage{babelbib}

\renewcommand{\baselinestretch}{1.0}

\begin{document}
\bibliographystyle{dcu}

\pagenumbering{arabic}

\begin{titlepage}
\begin{center}

{\textbf{ Model and heuristics for the Assembly Line Worker Integration and Balancing Problem}}

\vspace{0.5cm} \emph{\textbf{Mayron C\'esar de O. Moreira}}
\\  Instituto de Ci\^encias Matem\'aticas e de Computa\c{c}\~ao
\\  Universidade de S\~ao Paulo
\\
\vspace{0.5cm} \emph{\textbf{Crist\'obal Miralles}}
\\  ROGLE - Departamento Organizaci\'on de Empresas
\\  Universitat Polit\`ecnica de Val\`encia
\\
\vspace{0.5cm} \emph{\textbf{Alysson M. Costa}}
\\  Department of Mathematics and Statistics
\\  University of Melbourne

\end{center}

\vspace{2cm}

\noindent \textbf{Abstract} \\
{\small   }

We propose the Assembly Line Worker Integration and Balancing Problem (ALWIBP), a 
new assembly line balancing problem arising in lines with conventional and disabled workers. 
The goal of this problem is to maintain high productivity levels by minimizing the number of 
workstations needed to reach a given output, while integrating in the assembly line a number of 
disabled workers. Being able to efficiently manage a heterogeneous workforce is especially 
important in the current social context where companies are urged to integrate workers
with different profiles. In this paper we present mathematical models and heuristic methodologies that can help assembly 
line managers to cope with this additional complexity; demonstrating by means of a robust 
benchmark how this integration can be done with losses of productivity that are much lower 
than expected.

\noindent \textbf{Keywords:} Assembly line balancing; disabled workers; mathematical modeling.
\end{titlepage}

\renewcommand{\baselinestretch}{1.0}
\newpage

\section{Introduction}

According to the International Labour Organization (ILO), people with disabilities represent 
an estimated 10 percent of the world's population, including approximately 500 million 
of working age; being apparent that in the unemployment rates of the disabled 
are much higher than the average.

Employment is the main path for social inclusion and participation in modern societies. Having 
a job is not only the basis for the survival and stability for many individuals, but also 
a key way of accessing many rights as citizens. Therefore the welfare and the social inclusion 
of the disabled depend very much on the degree of labor integration they are able to achieve. 
Different active policies to fight against discrimination have been set during the last decades, 
following models that are more/less inclusive depending on the local culture.
Across specific national legislations, a general common formula is to reserve a share of
workplaces in ordinary companies for people with disabilities. This share normally increases 
with the size of the company and, depending on the country legislation, usually goes from 2\% 
to 5\% of the jobs.

Unfortunately, it is also a common phenomenon in many countries that this share is not 
always respected, indicating that the solution should come not only by legal 
imposition, but mainly by overcoming the prejudices about the capabilities of the disabled, 
and by the genuine commitment of ordinary companies to include integration programs in their 
strategies. The aim of this paper is to contribute in making this commitment easier: 1) 
by providing the production managers with practical approaches 
that ease the integration of disabled workers in the production lines; 2) by demonstrating 
that, through the approaches proposed, the productivity of production systems suffers 
little (and often none) decrease.

Once stated the great importance of integrating Disabled into the workforce of 
ordinary companies, we should make a brief introduction on some previous 
work inspired on the specific scenario of the so-called ``Sheltered Work Centers for 
Disabled'' (henceforth SWDs). SWDs are a special work formula legislated in many countries 
(with different variants) whose only difference from an ordinary company is that most of its workers 
must be disabled, and therefore they receive some institutional 
help in order to be able to compete in real markets. 
This labor integration formula has been successful in decreasing the former high unemployment rates of countries 
like Spain, and one of the strategies used by SWDs to facilitate the labor integration 
has been the adoption of assembly lines. In this sense, \cite{miralles07advantages} 
were the first to evidence how the integration of disabled workers in the productive 
systems can be done without losing, even gaining, productive efficiency through the 
use of assembly lines. This pioneer reference defined the so-called Assembly Line 
Worker Assignment and Balancing Problem (ALWABP) and demonstrated how the division
of worker into single tasks becomes a powerful tool for making certain workers
disabilities invisible.

\subsection{Contribution and outline of this work}

ALWABP was inspired in the SWDs reality, where the very high diversity of 
most of the workers and their limitations are the main characteristics. This 
scenario is quite different to that one of an ordinary company; where the aim 
is to efficiently integrate in the workforce just some workers, often to cope with the 
2\% to 5\% of disabled workers legislation requirements. In this case the problem 
supposes much less diversity in the input data, and can also be stated with very 
different approaches with respect to the ALWABP, regarding the objective function, 
the hypothesis and model defined, and the kind of appropriate solution procedures.

The aim of this paper is to introduce and analyze this 
new problem that has been named ``Assembly Line Worker Integration and Balancing Problem'' 
(ALWIBP). Our study aims to answer specific requirements that normally arise in assembly 
lines of ordinary companies, where only few disabled workers have to be integrated, providing 
the production managers with practical tools that ease the integration of disabled workers in 
the most efficient manner. We propose new mathematical models for the problem as well 
as a constructive heuristic based on the similarities between the proposed problem and the 
so-called Simple Assembly Line Balancing Problem (SALBP).

The remainder of this paper is structured as follows: in Section~\ref{ALWIBP}, we state 
a formal codification of the new problem and some extensions, analyzing their practical 
implications and reviewing references of the literature with useful related approaches. 
Section~\ref{sec:mathematical}  then presents the corresponding IP models for the proposed 
versions of the ALWIBP while Section~\ref{sec:heuristic} describes a fast heuristic that 
has been developed to solve the problem. A experimental study in order to analyze the 
effectiveness of the proposed models and algorithms is conducted in 
Section~\ref{sec:experimental}. General conclusions end this manuscript.

\section{The Assembly Line Worker Integration and Balancing Problem}
\label{ALWIBP}

\subsection{Introduction: SALBP vs ALWABP}

The SALBP was initially reviewed by \cite{baybars86survey} and consists of an assembly line balancing problem with several 
well-known simplifying hypotheses. This classical single-model problem which aims at
finding the best feasible assignment of tasks to stations so that certain precedence 
constraints are fulfilled, has been the reference problem in the literature in its two 
basic versions: when the cycle time $C$ is given, and the objective is to optimize the number 
of necessary workstations, the problem is called SALBP-1. Whereas when there is a given 
number $m$ of workstations, and the goal is to minimize the cycle time $C$ the literature 
knows this second version as SALBP-2 \citep{scholl99balancing}.

A trend in Assembly Line research in the last decade has been to narrow 
the gap between the theoretical proposals and the industrial reality,
which faces multiple specific configurations
such as multi-manned workstations \citep{Dimitriadis2006,kellegoz2012efficient},  
two sided assembly lines \citep{kim2009mathematical,ozcan2009multiple}, or operator allocation 
in job sharing and operator revisiting lines \citep{zeng2012operator}, among many others.
As part of this trend, \cite{miralles07advantages} properly defined the ALWABP, a generalization 
of the SALBP where, in addition to the assignment of tasks to 
stations, a set of heterogeneous workers also has to be assigned to stations. In this scenario 
each task has a worker-dependent processing time, which allows taking into account the limitations 
and specific production rates of each worker. Moreover, when the time to execute a task for certain worker is very 
high, this assignment is considered infeasible in the input data matrix.

Since \citep{miralles07advantages}, many other references have contributed 
to give ALWABP visibility throughout academia, proposing different 
methods to solve the problem. The same authors have later developed a 
branch-and-bound algorithm for the problem, obtaining the exact solution of 
small-sized instances \citep{miralles08branch}. Because of the problem 
complexity and the need to solve larger instances, the literature has since 
then shifted its efforts to heuristic methods.  The current state-of-the-art 
methods for solving the ALWABP are the  iterated beam search (IBS) metaheuristic 
of \cite{blum11solving}, the biased random-key genetic algorithm of 
\cite{moreira12simple}, the iterative genetic algorithm of \cite{mutlu13iterative},
the heuristic and the branch-and-bound algorithms of \cite{borba2013heuristic} and
the branch-and-bound algorithm of \cite{vila2014branch}. 

\subsection{ALWIBP}

The ALWABP problem was inspired in the SWDs reality with most workers presenting 
a high diversity of operation times; whereas the ALWIBP scenario introduced in Section 
1 intends to simulate the more inclusive situation in which 
disabled workers relative (in a small number) are integrated in a conventional assembly line. It has to be noted that
the main (and only studied) problem focusing on disabled integration 
in assembly lines has been the ALWABP-2 \citep[e.g.]{miralles08branch,moreira12simple}, since the 
typical objective at SWD is to be as efficient as possible with the (diverse) available workforce. 

In the scenario associated with the ALWIBP, it makes sense to deal with the type 1 problem, 
since a reasonable aim of a production manager can be to integrate the given disabled 
workers (in some cases some 2 or 5\% of workers, or even more whether some compensation is 
needed due to low shares in other factory sections) 
while minimizing the number of additional workstations needed for doing so. This problem is 
named ALWIBP-1, by analogy with the SALBP case.

In addition to this basic objective, once inside the solution subspace with minimal number 
of workstations, the manager may aim to find those assignments in which the idle time in 
stations with disabled workers is minimum, in order to increase their participation in the 
production process. We call this extension ALWIBP-1S$_{\textrm{min}}$. If, according to 
\cite{boysen07classification} classification, ALWABP-2 was stated as [$pa,link,cum|equip|c$], in this case
we can define ALWIBP-1 as [$pa,link,cum|equip|m$], while the ALWIBP-1S$_{\textrm{min}}$ can be stated as 
[$pa,link,cum|equip|m,{SLL}^{stat}$] using the same codification scheme.

In the following, we propose integer linear models for the basic ALWIBP-1 situation and also 
for the extension proposed. 

\section{Mathematical models}
\label{sec:mathematical}

In this section, we present a mathematical model for the ALWIBP-1 defined earlier, and 
further extend it to cope with ALWIBP1-$\textrm{S}_{\textrm{min}}$ extra objective with the 
use of the following notation: 

\vspace{0.3cm}
\noindent\begin{longtable}{lp{0.65\linewidth}}
  $N$ &set of tasks to be assigned; \\
  $S$ & set of workstations;\\
  $W$ & set of disabled workers, $|W| \le |S|$; \\
  $t_i$ & execution time of task $i$ when assigned to a ``conventional'' worker;\\
  $t_{wi}$ & execution time of task $i$ when assigned to disabled worker $w \in W$;\\
  $I_w \subseteq N$ & set of unfeasible tasks for worker $w \in W$;\\
  %$P_i = \{ j\mid (j,i)\in E\}$& set of immediate predecessors of task $i$;\\
  $F_i$ & set of immediate successors of task $i$;\\
  $F_i^*$ & set of all successors of task $i$.
\end{longtable}
\vspace{0.3cm}

\subsection{Model for ALWIBP-1}

The proposed formulation follows the idea used by \cite{petterson75assembly} 
when modeling the SALBP-1. Let $q$ be an artificial task and $D_q=\{i \in N| F_i = \varnothing \}$ 
be the set of tasks that have no followers. We define a new precedence graph in which $N' = N \cup \{q\}$
and all tasks in $D_q$ precede task $q$. The execution time of task $q$ is always 0, 
i.e., $t_q = t_{wq} = 0, \forall w\in W$. The ALWIBP-1 can thus be modelled as:

\begin{equation}
\label{ALWIBP1-1}
\textrm{Min }   \sum_{s\in S} sx_{sq}
\end{equation}

\hspace{1cm} subject to
\begin{alignat}{2}
\label{eq:ALWIBP1-2}
\sum_{s \in S} x_{si} = 1, &&  \qquad &\forall i \in N',\\
\label{eq:ALWIBP1-3}
\sum_{s \in S} y_{sw} = 1, && \qquad &\forall w \in W,\\
\label{eq:ALWIBP1-4}
\sum_{w\in W} y_{sw} \le 1, && \qquad &\forall s \in S, \\
\label{eq:ALWIBP1-5}
\sum_{s \in S| s \ge k} x_{si} \leq \sum_{s \in S| s \ge k} x_{sj}, && \qquad
&\forall i,j \in N', j \in F_i, k \in S, k \neq 1,\\
\label{eq:ALWIBP1-6}
\sum_{i\in N'} t_{i}x_{si} \leq \overline{C}, && \qquad &\forall s \in S,\\
\label{eq:ALWIBP1-7}
\sum_{i\in N'\backslash I_w} t_{wi}x_{si} \leq \overline{C} + L_w(1 - y_{sw}), && \qquad &\forall s \in S, \forall w \in W,\\
\label{eq:ALWIBP1-8}
y_{sw} \le 1 - x_{si}, && \qquad &\forall s \in S, \forall w \in W, \forall i
\in I_w,\\
\label{eq:ALWIBP1-9}
\sum_{s \in S| s \ge k}  y_{sw} \le \sum_{s \in S| s \ge k}  x_{sq}, &&  \qquad & \forall w \in W, \forall k \in S, k\neq 1,\\
\label{eq:ALWIBP1-10}
x_{si} \in \{0,1\}, && \quad &\forall s \in S, \forall i\in N',\\
\label{eq:ALWIBP1-11}
y_{sw} \in \{0,1\}, && \quad &\forall s \in S, \forall w\in W.
\end{alignat}

where:

$\begin{array}{ll}
x_{si}   & \textrm{binary variable equals to one if task } i \in N' \textrm{ is assigned to workstation } \\
	 & s \in S, \\
y_{sw}   & \textrm{binary variable equals to one if a disabled worker } w \in W \textrm{ is assigned to workstation } \\
	 & s \in S, \\
L_w   & \textrm{large constant}, w \in W.
\end{array}$
\vspace{0.3cm}

The objective function minimizes the index associated with the last station (the one 
that executes the fictitious last task $q$). 
In association with constraints (\ref{eq:ALWIBP1-3}) which state that all 
disabled workers are assigned, this objective function minimizes the number of ``conventional'' 
workers used in the line. Constraints (\ref{eq:ALWIBP1-4}) guarantee that each workstation
receives at most one (disabled) worker. This is a fair assumption in many practical
situations even for lines solely with ``conventional'' workers. In the case of workers with
disabilities, this fact gains in importance since these workers might have
special constraints which might force that the workstation be modified
accordingly. Constraints (\ref{eq:ALWIBP1-2}) ensure 
that all tasks are assigned, while
constraints (\ref{eq:ALWIBP1-5}) guarantee that the precedence relations are respected.
These inequalities were proposed by \cite{ritt2013improved} which analyzed several versions of precedence
constraints and concluded that constraints (\ref{eq:ALWIBP1-5}) presented better theoretical
and practical results. 

Constraints (\ref{eq:ALWIBP1-6}) and (\ref{eq:ALWIBP1-7}) ensure that the cycle time is respected 
at stations without and with disabled workers, respectively. Constant $L_w$ must be
sufficiently large to deactivate these last constraints if $y_{sw}=0$. Therefore, we take
$L_w=\sum_{i \in N\setminus I_w} |t_{wi}-t_i|$. This expression assumes the maximum additional
time that a disabled worker $w$ spends at a station in comparison to a conventional worker 
(this would be the additional time needed for the execution of all worker-feasible tasks).
Also, we highlight that even though a workstation must have only a single worker,
this set of constraints allow that more than one task can be assigned to the same station.

Finally, constraints (\ref{eq:ALWIBP1-8}) and (\ref{eq:ALWIBP1-9}) guarantee that disabled workers are not assigned 
to tasks  which they  are not able to execute and that they execute at least one task, respectively.

\subsection{ALWIBP-1S$_{\textrm{min}}$ model extension}

In order to extend this formulation to the ALWIBP-1S$_{\textrm{min}}$ case, one can simply note that this new problem is characterized by the addition of another term in the objective function related to the idle time of the disabled workers. The new goal is to hierarchically minimize the number of stations (with higher priority) and the idle time of the stations with disabled workers. Thereby, this version of the problem aims to obtain more balanced solutions
that increase the participation of these workers.

To model this situation, we use non-negative real variables $\delta_w, w \in W$, to measure the idle time of each disabled worker $w$. The values of these new 
variables are obtained with the aid of slack variables
$l_{sw}, \forall s \in S, \forall w \in W$ associated to constraints 
(\ref{eq:ALWIBP1-7}), which are rewritten as:
\begin{alignat}{2}
\label{eq:idle-1}
\sum_{i\in N'\backslash I_w} t_{wi}\cdot x_{si} + l_{sw} = \overline{C} + L_w(1 - y_{sw}), && \qquad &\forall s \in S, \forall w \in W,
\end{alignat}
\noindent and with  new constraints which are added to establish the correct relations between 
 $\delta_{w}$ and the slack variables:
\begin{alignat}{2}
\label{eq:idle-4}
\delta_w \geq l_{sw} - (1-y_{sw})\cdot\left(L_w + \sum_{i\in N'} t_i\right),  && \qquad &\forall s \in S, \forall w\in W.
\end{alignat}
\vspace{0.3cm}

Variables $l_{sw}$ measure the slack for constraints (\ref{eq:idle-1}) 
obtained with each pair $s\in S, w\in W$, while constraints (\ref{eq:idle-4}) define $\delta_w$ as the idle time 
of worker $w$, which is obtained with the value of $l_{sw}$ in the constraint corresponding to the 
actual assignment $y_{sw}$.
\vspace{0.2cm}

The objective function can now include terms associated with these idle times, and the ALWIBP-1S$_{\textrm{min}}$ can be written as:

\begin{equation}
\label{eq:ALWIBP-1Smin}
\textrm{Min }   \sum_{s \in S}  sx_{sq}   + \sum_{w\in W} \frac{\delta_{w}}{\overline{C}|W|}
\end{equation}
\hspace{1cm} subject to
\begin{alignat}{2}
(\ref{eq:ALWIBP1-2})-(\ref{eq:ALWIBP1-6}),  (\ref{eq:ALWIBP1-8})-(\ref{eq:idle-4}),\nonumber \\
\label{eq:ALWIBP1S-18}
l_{sw} \in \mathbb{R}_{+}, && \qquad &\forall s \in S, \forall w \in W, \\
\label{eq:ALWIBP1S-20}
\delta_{sw} \in \mathbb{R}_{+}, && \qquad &\forall s \in S, \forall w \in W.
\end{alignat}

The constant term multiplying the idle time variables imposes a hierarchical characteristic 
in the objective function, giving priority to the minimization of stations and using the idle 
times as a secondary objective.

\section{Constructive Insertion Heuristic}
\label{sec:heuristic}

We propose a Constructive Insertion Heuristic (CIH) for the ALWIBP-1. This heuristic  relies on the similarities between the ALWIBP and the classical SALBP. Indeed, during the heuristic procedures, a SALBP-1 solution is found and then iteratively adapted to incorporate the heterogeneous workers available in the ALWIBP-1. In the following section, this heuristic is presented in general terms. A formal description of the algorithm is presented in Section~\ref{sec:formal}. Finally, some variants of the procedure  and two post-optimization routines are described in Sections~\ref{sec:variants} and \ref{sec:post}, respectively. 

\subsection{General structure of the proposed heuristic}

The main steps of the proposed heuristic are described in Algorithm~\ref{alg:general}.

\begin{algorithm}[!h] \caption{Constructive insertion heuristic }
  \small
  \label{alg:general}
  \begin{algorithmic}[1]
  \STATE Ignore heterogeneous workers and obtain a SALBP-1 solution;
  \STATE Divide \textit{remaining} line in segments (containing a certain number of stations each);
  \STATE Try all available heterogeneous workers in each station of the first line segment;
  \STATE Select best assignment and fix the solution on stations prior to the selected one. 
  \STATE If there are still workers to assign, go to step 2; Otherwise, end.
  \end{algorithmic}
\end{algorithm}

%Fig.~\ref{fig:sketch} sketches this approach for a problem with two disabled workers. 
In the following subsections, each step of the algorithm is explained in details. 

\subsubsection{Step 1: obtain a SALBP-1 solution}
\label{step1}

A feasible solution for the SALBP-1, $x_{ref}$, is used as starting point of the algorithm. This solution can be obtained with any heuristic approach. In this study, we use the best solution found by CPLEX 12.4 in 10 minutes of running time on a single thread of a 2.3Ghz Intel machine with 6Gb of memory.

\subsubsection{Step 2: divide line into segments}
\label{step2}

The part of the line that has not yet been fixed is divided into $|W_a|$ segments, where $W_a$ is the set of disabled workers which have not yet been assigned. If a single worker is available, all stations in the remaining line are considered, otherwise, we consider stations in the current segment, $S_c$, defined as  $S_c \leftarrow\displaystyle\left\{s_b + 1,...,s_b + 1 + \left\lfloor\dfrac{m_{c} - s_{b}}{|W_a|}\right\rfloor\right\}$, where $s_b$ is the last fixed station and $m_c$ is the total number of stations in the current solution.

This choice of line segments follows two main rationales: the first is to have a sufficient number of stations for assigning the remaining workers with disabilities after each of these workers are assigned. The second is to somehow impose an even distribution of such workers along the line, with the goal of allowing a deeper integration between workers with and without disabilities - with the more experienced workers serving as ‘monitors’ for the newcomers.

\subsubsection{Step 3: Try all possible assignments in current segment}
\label{step3}

In this step, we try to obtain a solution by considering all assignments of workers 
in $W_a$ to stations in the current analyzed line segment $S_c$, one at a time. When a 
station $s\in S_c$ is considered for the assignment of worker $w\in W_a$, the following 
solution is obtained: in stations prior to $s$, the current solution, $x_{tmp}$, remains fixed. 
From station $s$ and on, the solution is rebuild using a constructive heuristic based on that of
\citep{scholl96simple}. We call this procedure \textit{salbp-1}($x^{tmp},w,s$). The 
following heuristic criteria were used to prioritize the assignment of tasks during the 
procedure (the reader is referred to the article of \cite{scholl96simple} for more details):

\begin{itemize}
 \item \textit{MaxTime}: Ascending order of task execution times, $t_i$;
 \item \textit{MaxPW}: Descending order of positional weights, ${pw}_{i} = t_i + \sum_{j \in F_i^{*}} t_j$;
  \item \textit{MaxIF}: Descending order of immediate followers, $F_i$;
  \item \textit{MaxF}: Descending order of followers, $F_i^{*}$;
\end{itemize}

In case a disabled worker is being considered, the appropriate parameters 
$t_{wi}$ and $t_{wj}$ are considered in criteria MaxTime and MaxPw.

\subsubsection{Step 4: Select best assignment}
\label{step4}

In this step, the best solution $x^*$ among those obtained for each assignment $(w,s)$, with $w\in W_a$  and $s\in S_c$ is chosen. A solution $x^1$ obtained with an assignment of worker $w_1$ to station $s_1$ is considered better than a solution $x^2$  ($x^{1} < x^{2}$) if the number of final stations of $x^1$ (designed by $stations(x^1)$) is smaller than $stations(x^2)$. If a tie happens, we consider the solution with larger idle time in the last station to be the best one.

%\begin{figure}[!h]
%\includegraphics[width=120mm]{sketch.pdf}
%\label{fig:sketch}
%\caption{Sketch of the steps of Algorithm~\ref{alg:general}}
%\end{figure}

\subsection{Formal structure of the proposed algorithm}
\label{sec:formal}

A formal description of the whole implemented heuristic is presented in Algorithm \ref{alg:formal}. It uses a given SALBP-1 solution as input as described in Section~\ref{step1}. Line division into segments described in Section~\ref{step2} is effected in lines 1-2 and in lines 20-21. The main loop of the algorithm occurs in lines 5-23: while there are workers yet to be assigned, each pair $(w\in W_a,s\in S_c)$ is tested by obtaining a complete solution with the procedure described in Section~\ref{step3} (line 11) and the best one (according to the criterion presented in Section~\ref{step4}) is selected (lines 12-15).

\begin{algorithm} \caption{Constructive Insertion Heuristic}
  \small
  \label{alg:formal}
  \begin{algorithmic}[1]
  \REQUIRE $x^{ref}$ (SALBP-1 reference solution);  
  \STATE $m_{c} \leftarrow$ \textit{stations}($x^{ref}$);
  \STATE $S_c \leftarrow \displaystyle\left\{1,...,\left\lceil\dfrac{m_{c}}{|W_a|}\right\rceil\right\}$;  
  \STATE $W_a \leftarrow W$;  
  \STATE Let $x^{*}$ be the incumbent solution (initially, $stations(x^{*}) = \infty$);
  \WHILE{$W_a \neq \varnothing$}
    \STATE $w_b \leftarrow \infty$;
    \STATE $s_b \leftarrow \infty$;
    \FORALL{$w \in W_a$}
      \FORALL{$s \in S_c$}
	\STATE Fix solutions of $x^{ref}$  on stations prior to $s$ in $x^{tmp}$;
	\STATE $x^c \leftarrow$ \textit{proc}\underline{ }\textit{salbp-1}($x^{tmp},w,s$);
	\IF{ $x^c < x^*$}
	  \STATE $w_b \leftarrow w$;
	  \STATE $s_b \leftarrow s$;
	  \STATE $x^* \leftarrow x^c$;
	\ENDIF
      \ENDFOR
    \ENDFOR
    %\IF{$s_b = \infty$ \textbf{ou} $w_b = \infty$}
     % \STATE Infeasible solution;
      %\STATE $W_a = \varnothing$;
    %\ELSE
      \STATE $W_a \leftarrow W_a \backslash \{w_b\}$;
      \STATE $m_{c} \leftarrow$ \textit{stations}($x^*$);
      \STATE $S_c \leftarrow \displaystyle\left\{s_b + 1,...,s_b + 1 + \left\lfloor\dfrac{m_{c} - s_{b}}{|W_a|}\right\rfloor\right\}$;  
      \STATE $x^{ref} \leftarrow x^{*}$;
    %\ENDIF
  \ENDWHILE
  \ENSURE $x^{*}$ (best solution found).
  \end{algorithmic}
\end{algorithm}

\subsection{Algorithmic variants}
\label{sec:variants}

The proposed algorithm tries to assign workers from the first ordered station to the last one. A variant of this strategy might start at the last stations and proceed backwards. In this case, the line division into segments presented in Section~\ref{step2} needs to be modified. Equations in line 2 and 21 of Algorithm~\ref{alg:formal} are modified for $S_c \leftarrow \displaystyle\left\{\left\lceil\dfrac{m_{c}}{|W_a|}\right\rceil,...,m_c\right\}$ and $S_c \leftarrow \displaystyle\left\{\left\lfloor\dfrac{s_{b} - 1}{|W_a|}\right\rfloor,..., s_b - 1\right\}$, respectively. Since in this backward approach the tasks and workers assigned in the end of the line are kept fixed, the algorithm may need to insert intermediate stations. Therefore, the tie-breaker for the procedure described in Section~\ref{step4} is no longer the idle time in the last station but the idle time in the last non-fixed station in the current solution.

In addition to changing the order the ordered stations are visited for worker assignments, one 
can also change the order the tasks are assigned in the procedure described in 
Section~\ref{step3}. Indeed, as usual, the precedence graph can be considered in its normal or reversed form, yielding four algorithmic variants (two worker assignment strategies combined with two task assignment strategies). Since the algorithm proposed is very fast, 
all four strategies are used at each run of the algorithm and the best result is kept; 
this increases the robustness of the algorithm with respect to the precedence
graph structure.

\subsection{Post-optimization routines useful for ALWIBP-1S$_{\textrm{min}}$}
\label{sec:post}

The solution obtained by the CIH is improved by means of a mixed-integer programming neighborhood based on the original model~(\ref{ALWIBP1-1})-(\ref{eq:ALWIBP1-11}). The main idea is to reduce the computational burden needed to solve the model by fixing the worker assignment variables and allowing changes only in task assignments variables. This can be done by solving the original model (\ref{ALWIBP1-1})-(\ref{eq:ALWIBP1-11}) with the addition of the following constraints:

\begin{equation}
y_{s(w),w} =1, \qquad \forall w \in W_a,
\end{equation}
where $s(w)$ indicates the station to which worker $w$ is assigned in the solution of the CIH.  Clearly, the model can be adjusted accordingly, to eliminate useless constraints and variables.

In order to further reduce the computational effort needed to solve the model, one can reduce the freedom of task assignment variables in the MIP neighborhood, by allowing tasks to be assigned only to the same station it was assigned in the CIH solution or to neighboring stations. This can be easily done with the addition of the following constraints:

\begin{alignat}{2}
\label{posOtm1_12}
x_{s(i),i} + x_{s(i)-1,i} + x_{s(i)+1,i} = 1, && \qquad &\forall i \in N, s(i) \in S\backslash \{1,m\},\\
\label{posOtm1_13}
x_{1i} + x_{2i} = 1, && \qquad &\forall i \in N, s(i)=1,\\
\label{posOtm1_14}
x_{mi} + x_{m-1,i} = 1, && \qquad &\forall i \in N, s(i) = m.
\end{alignat}
Where $s(i)$ indicates the station to which task $i$ was assigned in the CIH solution.

Additionally, the objective function of this post-optimization program can be changed to incorporate the characteristics of the ALWABP-1$\textrm{S}_{\textrm{min}}$, in which one prioritizes solutions with low idle times for the disabled workers.  In this case, it suffices to use  (\ref{eq:ALWIBP-1Smin}) as objective function and add constraints:

\begin{equation}
\sum_{i\in N} p_{s(w),i}\cdot x_{si} + \delta_w = C,  \qquad \forall w\in W,\\
\end{equation}
Note that these simpler constraints replace constraints (\ref{eq:idle-1})-(\ref{eq:idle-4}) in the 
original ALWIBP-1$\textrm{S}_{\textrm{min}}$ problem, since the worker assignments are 
already known in the post-optimization phase.

\section{Experimental study}\label{sec:experimental}

\subsection{Justification of a new ALWIBP benchmark}

As discussed in section 2, the ALWABP was inspired by the situation found at SWDs where the 
very high diversity of workers and their limitations are the main characteristics; whereas 
the ALWIBP scenario pretends to simulate the ``desirable'' situation of 
only a small percentage of disabled workers being integrated in conventional assembly lines. Moreover, 
as stated earlier, the main (and only studied) approach in this scenario has been ALWABP-2, since 
the typical objective in SWD is to be as efficient as possible with the (diverse) available 
workforce. Instead, in this new scenario, the assembly line balancing of type 1 (minimization of the number
of stations given the desired cycle time) approach becomes realistic, since the basic 
aim of a production manager can be to integrate the normative (common in most countries) 2\%-5\% 
of disabled workers into the assembly line; while 
maintaining a given productivity. 
Many previous proposals for the ALWABP-2 were evaluated with the set of 320 benchmark instances 
first proposed by \cite{chaves09hybridb}. Once stated the completely different scenario where 
the ALWIBP arises, it is clear that this classical ALWABP benchmark is not useful here since, as 
explained above: (1) only a little share of the workers are disabled; and (2) the basic aim 
is now to minimize the number of workstations with non-disabled workers (ALWIBP-1 perspective). 

\subsection{ALWIBP benchmark scheme}
As many other ALB approaches, the ALWABP benchmark was constructed from the only SALBP 
reference (the Scholl data collection of www.assembly-line-balancing.de); 
that was considered robust enough and has been extensively used to test most
 proposals in the literature so far. Nevertheless, as recently demonstrated 
by \cite{otto13salpbgen}, this framework does not seem rigorous enough. The problems 
were collected from different empirical and not empirical sources, and are based  
on only 25 precedence graphs; where just 18 distinct graphs have more than 25 tasks and thus 
are meaningful for comparing solution methods. Moreover, \cite{otto13salpbgen} also point out the 
triviality of many classical instances. Therefore, they propose a SALBP generator and a new very robust challenging benchmark whose graphs 
morphologies include a sufficient variety of chains, bottlenecks and modules. Basically, it 
has different cells of data sets (with 25 different instances per cell) following a 
full-factorial design for the following parameters:

% , proposing to
% avoid these two inconvenient characteristics (which is: low diversity of graphs structure, and
% triviality), a SALBP generator and a new very robust challenging benchmark whose graphs 
% morphologies include a sufficient variety of chains, bottlenecks and modules. Basically, they 
% propose different cells of data sets (with 25 different instances per cell) following a 
% full-factorial design for the following parameters:

\begin{enumerate}
\item Type of the graph: precedence graphs containing more chains, more bottlenecks, or
a mix of both.
\item Order Strength: ``low'', ``medium'' and ``high''.
\item Distribution of task times: ``peak at the bottom'', ``bimodal'' and ``peak in the middle''.
\item Number of tasks: ``small'', ``medium'', ``large'' and ``very large''.
\end{enumerate}

Thus, we selected as basis the following collection of subsets from the
\cite{otto13salpbgen} benchmark: 

\begin{enumerate}
\item We consider that diversity of graphs is sufficiently ensured selecting only the ``mixed''
instances (that have both chains and bottlenecks).
\item From them, we only select the data subsets with ``low'' and ``high'' Order Strength (that
would give us clearer correlations if needed).
\item And from them, regarding distribution of task times, we select the subsets with ``peak
at the bottom'' and ``bimodal'' distribution (we discard the ``peak in the middle'' subset
because the optimal SALBP solution is unknown for many of the instances).
\item Finally, we take the resulting 100 robust instances of the ``medium'' (with $n = 50$ tasks),
``large'' (with $n = 100$ tasks), and ``very large'' (with $n = 1000$ tasks) subsets, and we
generate three corresponding benchmarks that we use separately in our experimentation
(we discard the ``small'' subset since it is advised for simple testing only).
\end{enumerate}

In all three cases we generate the benchmark using the same procedure: for each original instance we respect 
the precedence network and the conventional task time, and then we generate four different 
instances by adding one disabled worker with: high or low variability of task time respect 
to the original ones, and high or low percentage of incompatibilities. The two levels defined 
for the task times variability used the distributions $U[t_i, 2t_i]$ and $U[t_i, 5t_i]$ for low 
and high variability, and the low and high percentage of incompatibilities in the tasks-workers 
matrix was set to 10\% and 20\% approximately. Following the same scheme we created 400 
additional instances with two workers, then three workers, and finally four workers. 

Following this outline, we finally obtain three reliable benchmarks with the structure described, one with 
1600 ``medium'' instances, a second one with 1600 ``large'' instances, and a third one with
1600 ``very large'' instances. It is important to note that in our source benchmark the 
instances are classified from ``less tricky'' to ``extremely tricky'' and it happens that, 
the four cells finally selected as base have a quite symmetric composition regarding this ``triviality'' measure.

Finally, it is also important to note that the following ALWIBP computational study has always
used an input cycle time of 1000 time units simply because it is the value used by 
\cite{otto13salpbgen}, making the results of both cases comparable.

\subsection{ALWIBP computational study}

% One basic aim of every company should be to integrate at least the normative percentage of 
% disabled workers into the workforce. In this experimental study, we aim to demonstrate that 
% the proposed methods enable the inclusion of higher percentages of disabled workers in 
% the line without important losses in productivity. Thus, productivity always means 
% somehow (output result / input resources involved) and in this case productivity can
% be defined as (throughput rate / number of workstations). As in all three benchmarks the 
% throughput rate is fixed (because cycle time is always 1000 time units), an increase 
% on the number of workstations means, in general, a decrease of productivity.
% 
% Thus, in a first experiment we use the ``medium'' (with $n=50$ tasks) and the ``large''data subset 
% (with $n=100$ tasks); obtaining exact solutions for almost all the instances 
% with the model(and extension) of section 3. In a second experiment they are compared 
% with the heuristic results. And in a third experiment we use the ``very large'' data
% subset (with $n=1000$ tasks) to get further evidence of the good behavior of the heuristic against 
% big problems that are difficult to solve exactly.

Our computational study consists of three parts: in the first experiment we use the ``medium'' and
``large'' benchmarks, obtaining exact solutions for almost all the instances with the model (and
extension) of section 3. In the second experiment we use the same two benchmarks to compare
the exact results with the ones obtained by the heuristic procedures. In the third experiment we
use the ``very large'' benchmark to get further evidence of the good behavior of the heuristic
against large problems that are difficult to solve exactly. Furthermore, we validate 
the obtained results by comparing them to those of another heuristic procedure.

\subsubsection{Experiment 1: exact results for ALWIBP-1 and ALWIBP-1S$_{min}$}
\label{sec:experiment1}

% As an initial experiment to check what would be the shape of this expected loss of  
% productivity, we take the solutions of the 100 selected SALBP ``medium'' ($n=50$) and
% ``large'' ($n=100$) instances from \cite{otto13salpbgen} that are the base of our
% benchmarks, and compare them with our solutions once 
% integrated (applying the model and extension of section 3) consecutively one, 
% two, three and four workers. 

One basic aim of every company should be to integrate at least the normative percentage of
disabled workers into the workforce. In this first experiment, we aim to demonstrate that the
proposed methods enable the inclusion of higher percentages of disabled workers in the line
without important losses in productivity. It has to be noted that productivity always means
somehow (output result / input resources involved), and in assembly lines productivity can
be defined as (throughput rate / number of workstations). Therefore (considering that in this
experiment the Throughput rate is fixed), an increase on the number of workstations means, in
general, a decrease of productivity.

Having this in mind, this first experiment was devoted to check what would be the shape of this
expected loss of productivity. For the medium and large instances, both ALWIBP-1 and ALWIBP-1S$_{min}$ models could be solved in 
reasonable computation times using the commercial 
package CPLEX 12.4 (1 \textit{thread}, limit time of 1800 sec, 
and 6Gb as tree limit) with Intel Core i7 3.4 GHz, 16Gb RAM. 

In Tables \ref{tab:models_alwibp_middle_inst} and \ref{tab:models_alwibp_large_inst}, the 
columns indicate:

\begin{itemize}
\item $\mathbf{\Delta}$: number of instances solved to optimality;
\item \textbf{t(s)}: computational time (on average);
\item $\mathbf{m_{\uparrow}}$ and $\mathbf{\sigma_{m_{\uparrow}}}$: number of additional workstations needed in the ALWIBP-1 solution, with respect to the best know solution of SALBP-1 (average and deviation);
\item $\mathbf{m_{\uparrow}}$ and $\mathbf{\sigma_{m_{\uparrow}}}(\%)$: percentage of additional workstations needed in the ALWIBP-1 solution, with respect to the best know solution of SALBP-1 (average and deviation);
\item $\mathbf{{\tau}}$: \textit{idle time} of stations with disabled workers (on average) in ALWIBP-1;
\item $\mathbf{{\tau}_{Smin}}$: \textit{idle time} of stations with disabled workers (on average) in ALWIBP-1S$_{min}$;
\item $\mathbf{\eta} (\%)$: percentage of tasks performed by disabled workers with respect to the mean number of tasks assigned to ordinary workers in ALWIBP-1;
\item $\mathbf{{\eta}_{Smin}} (\%)$: percentage of tasks performed by disabled workers with respect to the mean number of tasks assigned to ordinary workers in ALWIBP-1S$_{min}$.
\item $\mathbf{\beta}$ and $\mathbf{\sigma_{\beta}} (\%)$: percentage of the number of disabled in the assembly line $\displaystyle\left(\dfrac{|W|}{m}\right)$, where $m$ is the number of stations in the ALWIBP-1 solution (average and deviation);
\item $\mathbf{\theta}$: number of instances which there is no increase of the number of stations with respect to the SALBP-1 solution;
%\item $\mathbf{\beta_{\theta}}$ and $\mathbf{\sigma_{\beta_{\theta}}}(\%)$: criterion $\mathbf{\beta}$ computed for instances in which there is 
no increase in the number of stations with respect to the SALBP-1 solution (average and deviation).
\end{itemize}

% Tabela da comparacao dos modelos ALWIBP-1 e ALWIBP-1Smin com instancias de 50 tarefas  %%%%%%%%%%%%%%%%%
\begin{sidewaystable}[htbp]
  \scriptsize
  \centering
  \caption{Computational results: ALWIBP-1 and ALWIBP-1S$_{min}$ models ($|N|$ = 50)}
    \begin{tabular}{ccc|cc|cc|cc|cc|c|c}
    \toprule
\textbf{|W|} & \textbf{Var} & \textbf{Inc} & $\mathbf{\Delta}$ & \textbf{t(s)} & $\mathbf{m_{\uparrow}} \pm \mathbf{\sigma_{m_{\uparrow}}}$ & $\mathbf{m_{\uparrow}} \pm \mathbf{\sigma_{m_{\uparrow}}} (\%)$ & $\mathbf{\tau}$ & $\mathbf{{\tau}_{S_{min}}}$ & $\mathbf{\eta} (\%)$ & $\mathbf{{\eta}_{S_{min}}} (\%)$ & $\mathbf{\beta}  \pm \mathbf{\sigma_{\beta}} (\%)$ & $\mathbf{\theta}$ \\    \midrule
    \multirow{3}[4]{*}{\textbf{1}} & \multirow{2}[2]{*}{\textbf{U[t,2t]}} & \textbf{10\%} & 95    & 104.0 & 0.2 $\pm$ 0.4 & 2.8\% $\pm$ 5.5\% & 149.9 & 1.6   & 66.2\% & 87.6\% & 10.9\% $\pm$ 2.7\% & 78 \\
    \textbf{} & \textbf{} & \textbf{20\%} & 95    & 90.6  & 0.2 $\pm$ 0.4 & 2.9\% $\pm$ 5.6\% & 131.4 & 3.1   & 62.8\% & 85.6\% & 10.9\% $\pm$ 2.7\% & 78  \\
    \textbf{} & \multirow{2}[2]{*}{\textbf{U[t,5t]}} & \textbf{10\%} & 95    & 125.0 & 0.4 $\pm$ 0.5 & 5.0\% $\pm$ 6.3\% & 231.6 & 9.6   & 43.7\% & 63.4\% & 10.6\% $\pm$ 2.6\% & 59  \\
    \textbf{} & \textbf{} & \textbf{20\%} & 95    & 120.4 & 0.4 $\pm$ 0.5 & 5.4\% $\pm$ 6.5\% & 244.1 & 10.2  & 39.9\% & 58.9\% & 10.6\% $\pm$ 2.5\% & 57  \\ \hline
    \multirow{3}[4]{*}{\textbf{2}} & \multirow{2}[2]{*}{\textbf{U[t,2t]}} & \textbf{10\%} & 92    & 177.7 & 0.4 $\pm$ 0.5 & 4.7\% $\pm$ 6.5\% & 105.9 & 1.9   & 67.9\% & 81.6\% & 21.3\% $\pm$ 5.0\% & 63  \\
    \textbf{} & \textbf{} & \textbf{20\%} & 96    & 85.8  & 0.4 $\pm$ 0.5 & 4.7\% $\pm$ 6.3\% & 101.3 & 3.5   & 63.0\% & 79.9\% & 21.3\% $\pm$ 5.2\% & 62  \\
    \textbf{} & \multirow{2}[2]{*}{\textbf{U[t,5t]}} & \textbf{10\%} & 93    & 166.3 & 0.7 $\pm$ 0.5 & 8.5\% $\pm$ 6.8\% & 127.1 & 7.4   & 51.6\% & 64.3\% & 20.5\% $\pm$ 4.8\% & 30  \\
    \textbf{} & \textbf{} & \textbf{20\%} & 93    & 152.7 & 0.7 $\pm$ 0.5 & 8.6\% $\pm$ 7.1\% & 140.1 & 8.1   & 47.7\% & 60.6\% & 20.5\% $\pm$ 4.8\% & 30  \\ \hline
     \multirow{3}[4]{*}{\textbf{3}} & \multirow{2}[2]{*}{\textbf{U[t,2t]}} & \textbf{10\%} & 91    & 195.7 & 0.5 $\pm$ 0.5 & 6.4\% $\pm$ 6.5\% & 72.8  & 2.7   & 73.0\% & 82.0\% & 31.4\% $\pm$ 7.3\% & 48  \\
    \textbf{} & \textbf{} & \textbf{20\%} & 89    & 213.6 & 0.6 $\pm$ 0.5 & 6.9\% $\pm$ 6.4\% & 76.5  & 3.1   & 69.0\% & 79.2\% & 31.3\% $\pm$ 7.3\% & 43  \\
    \textbf{} & \multirow{2}[2]{*}{\textbf{U[t,5t]}} & \textbf{10\%} & 95    & 109.6 & 1.1 $\pm$ 0.5 & 13.3\% $\pm$ 8.1\% & 114.0 & 7.4   & 51.3\% & 62.0\% & 29.4\% $\pm$ 6.5\% & 8  \\
    \textbf{} & \textbf{} & \textbf{20\%} & 90    & 217.0 & 1.2 $\pm$ 0.5 & 14.1\% $\pm$ 8.6\% & 126.5 & 11.0  & 51.1\% & 57.9\% & 29.2\% $\pm$ 6.3\% & 6  \\ \hline
    \multirow{3}[4]{*}{\textbf{4}} & \multirow{2}[2]{*}{\textbf{U[t,2t]}} & \textbf{10\%} & 89    & 222.7 & 0.7 $\pm$ 0.5 & 8.3\% $\pm$ 6.8\% & 56.5  & 3.9   & 75.9\% & 81.3\% & 41.2\% $\pm$ 9.5\% & 32 \\
    \textbf{} & \textbf{} & \textbf{20\%} & 86    & 301.2 & 0.8 $\pm$ 0.5 & 9.2\% $\pm$ 7.5\% & 68.0  & 5.0   & 70.5\% & 77.7\% & 40.8\% $\pm$ 9.3\% & 28 \\
    \textbf{} & \multirow{2}[2]{*}{\textbf{U[t,5t]}} & \textbf{10\%} & 89    & 278.8 & 1.5 $\pm$ 0.6 & 17.0\% $\pm$ 9.7\% & 109.6 & 11.7  & 53.7\% & 59.7\% & 37.9\% $\pm$ 7.8\% & 1 \\
    \textbf{} & \textbf{} & \textbf{20\%} & 84    & 362.9 & 1.6 $\pm$ 0.6 & 18.1\% $\pm$ 9.9\% & 120.5 & 11.7  & 50.7\% & 58.2\% & 37.6\% $\pm$ 7.9\% & 0  \\
    \bottomrule
    \end{tabular}%
  \label{tab:models_alwibp_middle_inst}%
\end{sidewaystable}%%%%%%%%%%%%%%%%%%%%%%%%%%%%%%%%%%%%%%%%%%%%%%%%%

% Tabela da comparacao dos modelos ALWIBP-1 e ALWIBP-1Smin com instancias de 100 tarefas %%%%%%%%%%%%%%
\begin{sidewaystable}[htbp]
  \scriptsize
  \centering
  \caption{Computational results: ALWIBP-1 and ALWIBP-1S$_{min}$ models ($|N|$ = 100)}
    \begin{tabular}{ccc|cc|cc|cc|cc|c|c}
    \toprule
    \textbf{|W|} & \textbf{Var} & \textbf{Inc} & $\mathbf{\Delta}$ & \textbf{t(s)} & $\mathbf{m_{\uparrow}} \pm \mathbf{\sigma_{m_{\uparrow}}}$ & $\mathbf{m_{\uparrow}} \pm \mathbf{\sigma_{m_{\uparrow}}} (\%)$ & $\mathbf{\tau}$ & $\mathbf{{\tau}_{S_{min}}}$ & $\mathbf{\eta} (\%)$ & $\mathbf{{\eta}_{S_{min}}} (\%)$ & $\mathbf{\beta}  \pm \mathbf{\sigma_{\beta}} (\%)$ & $\mathbf{\theta}$ \\
    \midrule
    \multirow{3}[4]{*}{\textbf{1}} & \multirow{2}[2]{*}{\textbf{U[t,2t]}} & \textbf{10\%} & 86    & 309.1 & 0.2 $\pm$ 0.4 & 1.0\% $\pm$ 2.3\% & 68.9  & 1.1   & 74.4\% & 87.0\% & 5.6\% $\pm$ 1.4\% & 83 \\
    \textbf{} & \textbf{} & \textbf{20\%} & 85    & 330.6 & 0.2 $\pm$ 0.4 & 1.0\% $\pm$ 2.3\% & 85.6  & 1.7   & 67.0\% & 85.7\% & 5.6\% $\pm$ 1.4\% & 82  \\
     \textbf{} & \multirow{2}[2]{*}{\textbf{U[t,5t]}} & \textbf{10\%} & 81    & 435.2 & 0.3 $\pm$ 0.5 & 2.0\% $\pm$ 2.9\% & 109.1 & 3.4   & 56.8\% & 66.0\% & 5.5\% $\pm$ 1.3\% & 66  \\
    \textbf{} & \textbf{} & \textbf{20\%} & 74    & 521.9 & 0.4 $\pm$ 0.5 & 2.4\% $\pm$ 3.1\% & 125.1 & 5.3   & 49.7\% & 67.0\% & 5.5\% $\pm$ 1.3\% & 59  \\ \hline
    \multirow{3}[4]{*}{\textbf{2}} & \multirow{2}[2]{*}{\textbf{U[t,2t]}} & \textbf{10\%} & 67    & 675.8 & 0.4 $\pm$ 0.5 & 2.3\% $\pm$ 3.1\% & 70.4  & 1.5   & 71.4\% & 88.5\% & 11.0\% $\pm$ 2.7\% & 61  \\
    \textbf{} & \textbf{} & \textbf{20\%} & 65    & 679.5 & 0.4 $\pm$ 0.5 & 2.4\% $\pm$ 3.1\% & 72.7  & 1.7   & 71.3\% & 85.5\% & 11.0\% $\pm$ 2.7\% & 59  \\
    \textbf{} & \multirow{2}[2]{*}{\textbf{U[t,5t]}} & \textbf{10\%} & 53    & 971.1 & 0.7 $\pm$ 0.5 & 4.2\% $\pm$ 2.9\% & 95.9  & 5.5   & 57.1\% & 65.9\% & 10.8\% $\pm$ 2.6\% & 28  \\
    \textbf{} & \textbf{} & \textbf{20\%} & 49    & 1045.3 & 0.8 $\pm$ 0.4 & 4.5\% $\pm$ 2.7\% & 105.5 & 4.8   & 55.4\% & 66.2\% & 10.8\% $\pm$ 2.6\% & 21  \\ \hline
    \multirow{3}[4]{*}{\textbf{3}} & \multirow{2}[2]{*}{\textbf{U[t,2t]}} & \textbf{10\%} & 47    & 1070.7 & 0.6 $\pm$ 0.5 & 3.3\% $\pm$ 3.1\% & 48.8  & 1.7   & 74.6\% & 85.2\% & 16.4\% $\pm$ 3.9\% & 43  \\
    \textbf{} & \textbf{} & \textbf{20\%} & 47    & 1099.7 & 0.6 $\pm$ 0.5 & 3.5\% $\pm$ 3.0\% & 55.0  & 2.8   & 73.4\% & 83.8\% & 16.4\% $\pm$ 3.9\% & 39  \\
    \textbf{} & \multirow{2}[2]{*}{\textbf{U[t,5t]}} & \textbf{10\%} & 28    & 1397.5 & 1.1 $\pm$ 0.5 & 6.4\% $\pm$ 3.2\% & 86.7  & 5.5   & 57.4\% & 65.1\% & 15.9\% $\pm$ 3.7\% & 5  \\
    \textbf{} & \textbf{} & \textbf{20\%} & 28    & 1407.1 & 1.1 $\pm$ 0.5 & 6.5\% $\pm$ 3.2\% & 87.4  & 4.1   & 58.0\% & 63.3\% & 15.9\% $\pm$ 3.7\% & 4  \\ \hline
    \multirow{3}[4]{*}{\textbf{4}} & \multirow{2}[2]{*}{\textbf{U[t,2t]}} & \textbf{10\%} & 34    & 1300.6 & 0.8 $\pm$ 0.5 & 4.4\% $\pm$ 2.9\% & 55.2  & 2.7   & 74.1\% & 82.6\% & 21.6\% $\pm$ 5.1\% & 26  \\
    \textbf{} & \textbf{} & \textbf{20\%} & 32    & 1362.7 & 0.8 $\pm$ 0.5 & 4.6\% $\pm$ 2.9\% & 55.5  & 3.4   & 72.4\% & 80.8\% & 21.6\% $\pm$ 5.1\% & 22  \\
    \textbf{} & \multirow{2}[2]{*}{\textbf{U[t,5t]}} & \textbf{10\%} & 8     & 1707.1 & 1.5 $\pm$ 0.5 & 8.4\% $\pm$ 4.0\% & 89.3  & 7.4   & 57.7\% & 61.4\% & 20.8\% $\pm$ 4.7\% & 2  \\
    \textbf{} & \textbf{} & \textbf{20\%} & 12    & 1724.5 & 1.5 $\pm$ 0.5 & 8.5\% $\pm$ 3.7\% & 109.6 & 8.2   & 57.6\% & 58.5\% & 20.8\% $\pm$ 4.8\% & 1  \\
    \bottomrule
    \end{tabular}%
  \label{tab:models_alwibp_large_inst}%
\end{sidewaystable}%%%%%%%%%%%%%%%%%%%%%%%%%%%%%%%%%%%%%%%%%%%%%%%%

We analyze first the results for ``middle'' instances. 
As expected, the increase in the number of stations grows with the number of disabled workers 
to be integrated and with the variability of the task times. Nevertheless, 
it can be observed that even in the most constrained case (4 disabled 
workers with execution times of up to 5 times the conventional time and 
20\% incompatibility), an average of only 1.6 new stations had to be added 
to integrate the workers.

This result is even more remarkable when the size of lines is taken into consideration. 
Indeed, in column 12, the $\mathbf{\beta}$ metric shows us that 
an average increase of 0.4 workstations imply in approximately
20\% of the workforce composed by disabled workers. This is, in fact,
relevant in real contexts, once the legislation of many countries requires a rate of
disabled workers ranging from 2\% to 5\% in industries.

Concerning the $\mathbf{\theta}$ metric, we see that in 38\% of the cases we obtain a line
balancing for the ALWIBP-1 case with the same number of workstations of the SALBP-1 solution, resulting
in the best of possible scenarios (integrate all disabled workers without increasing
production costs). Furthermore, the percentage of disabled workers that take part in the
workforce keep at least in 5\% in these cases.

As the sub-space of possible ALWIBP-1 optimal solutions is large, we can combine this primary aim of minimizing 
conventional workstations with the (important) secondary objective of minimizing the 
idle time of disabled workers. After applying the extension ALWIBP-1S$_{min}$ to all 
the instances the results show that we maintain in all cases the same rough productivity 
(no increase of the number of workstations), while we reach a very little mean idle 
time (see column $\mathbf{{\tau}_{S_{min}}}$).
We note that in the column related to the $\mathbf{\eta} (\%)$ and 
$\mathbf{{\eta}_{Smin}} (\%)$ criteria, since the
increase of these values follows the reduction of idle time in all scenarios studied.

Concerning the ALWIBP models in the ``large'' instances (with 100 tasks), Table 
\ref{tab:models_alwibp_large_inst} has shown
that the obtention of optimal solutions becomes difficult when 3 and 4 disabled workers
must be inserted in assembly lines (column $\Delta$). However, the number of additional
workstations is smaller than we notice in the experiments of ``middle'' instances, which
implies that its approach may be useful in more complex assembly lines. Furthermore, the percentage
of disabled workforce integrated continues in at least 5\%, as it is desired.

\subsubsection{Experiment 2: Heuristics comparison}

For validation purposes, we solved the same 1600 ``medium'' and ``large'' instances with the 
CIH procedure and its variations with post-optimization phase (CIH-LS1 and CIH-LS2) 
presented in section 4. The MIP-local search routines were run with the package CPLEX
12.4 with the same parameter settings, except for the computational time limit, which
was set to 60s.

In Tables \ref{tab:heuristic_cih_small_inst} and \ref{tab:heuristic_cih_large_inst}, the columns indicate:

\begin{itemize}

\item $\mathbf{m_{h_{\uparrow}}}$ and $\mathbf{\sigma_{m_{h_{\uparrow}}}}$: number of additional workstations in the CIH solution, with respect to the solution obtained when solving the ALWIBP-1 model (average and deviation);
\item $\mathbf{m_{h_{\uparrow}}}$ and $\mathbf{\sigma_{m_{h_{\uparrow}}}}(\%)$: percentage of additional workstations in the CIH solution, with respect to the solution obtained when solving the ALWIBP-1 model (average and deviation);
\item $\mathbf{t_{h}(s)}$: Computational time used by the heuristic (on average);
\item \textbf{Ties}: Number of instances in which the CIH solution had the same number of workstations (no increase) as the ALWIBP-1 model solution;
\item \textbf{Ties}$_{sp1}$: Number of instances in which the CIH solution had the same number of workstations (no increase) as the SALBP-1 model solution;
\end{itemize}

% Table generated by Excel2LaTeX from sheet 'Heuristicas_ALWIBP_Resumo' %%%%%%%%%%%%%%%%%%%%%%%%%%
\begin{sidewaystable}[htbp]
  \scriptsize
  \centering
    \caption{Computational results: CIH approaches ($|N|$ = 50)}
    \begin{tabular}{ccc|cccc|cccc|cccc}
    \toprule
    \multirow{2}[0]{*}{\textbf{|W|}} & \multirow{2}[0]{*}{\textbf{Var}} & \multirow{2}[0]{*}{\textbf{Inc}} & \multicolumn{4}{c}{\textbf{CIH}} & \multicolumn{4}{c}{\textbf{CIH-LS1}} & \multicolumn{4}{c}{\textbf{CIH-LS2}} \\
    \cline{4-15}
          &       &       & $\mathbf{m_{h_{\uparrow}}} \pm \mathbf{\sigma_{m_{h_{\uparrow}}}}$ & $\mathbf{m_{h_{\uparrow}}} \pm \mathbf{\sigma_{m_{h_{\uparrow}}}} (\%)$ & $\mathbf{t_{h}(s)}$ & \textbf{Ties} & $\mathbf{m_{h_{\uparrow}}} \pm \mathbf{\sigma_{m_{h_{\uparrow}}}}$ & $\mathbf{m_{h_{\uparrow}}} \pm \mathbf{\sigma_{m_{h_{\uparrow}}}} (\%)$ & $\mathbf{t_{h}(s)}$ & \textbf{Ties} & $\mathbf{m_{h_{\uparrow}}} \pm \mathbf{\sigma_{m_{h_{\uparrow}}}}$ & $\mathbf{m_{h_{\uparrow}}} \pm \mathbf{\sigma_{m_{h_{\uparrow}}}} (\%)$ & $\mathbf{t_{h}(s)}$ & \textbf{Ties} \\ \hline
    \multirow{3}[4]{*}{\textbf{1}} & \multirow{2}[2]{*}{\textbf{U[t,2t]}} & \textbf{10\%} & 0.2 $\pm$ 0.4 & 1.8\% $\pm$ 4.3\% & 0.0     & 84    & 0.1 $\pm$ 0.2 & 0.6\% $\pm$ 2.7\% & 10.5  & 95    & 0.1 $\pm$ 0.3 & 1.3\% $\pm$ 3.6\% & 0.1   & 88 \\
    \textbf{} & \textbf{} & \textbf{20\%} & 0.1 $\pm$ 0.3 & 1.3\% $\pm$ 3.6\% & 0.0     & 88    & 0.0 $\pm$ 0.1 & 0.2\% $\pm$ 1.7\% & 10.8  & 98    & 0.1 $\pm$ 0.3 & 1.1\% $\pm$ 3.4\% & 0.1   & 89 \\
    \textbf{} & \multirow{2}[2]{*}{\textbf{U[t,5t]}} & \textbf{10\%} & 0.2 $\pm$ 0.4 & 2.1\% $\pm$ 4.5\% & 0.0     & 81    & 0.1 $\pm$ 0.3 & 0.6\% $\pm$ 2.6\% & 10.8  & 92    & 0.2 $\pm$ 0.4 & 1.8\% $\pm$ 4.2\% & 0.1   & 83 \\
    \textbf{} & \textbf{} & \textbf{20\%} & 0.2 $\pm$ 0.4 & 2.3\% $\pm$ 4.7\% & 0.0     & 79    & 0.0 $\pm$ 0.2 & 0.4\% $\pm$ 1.8\% & 9.9   & 96    & 0.2 $\pm$ 0.4 & 2.0\% $\pm$ 4.4\% & 0.0     & 81 \\ \hline
    \multirow{3}[4]{*}{\textbf{2}} & \multirow{2}[2]{*}{\textbf{U[t,2t]}} & \textbf{10\%} & 0.2 $\pm$ 0.4 & 2.3\% $\pm$ 4.6\% & 0.0     & 79    & 0.1 $\pm$ 0.3 & 0.9\% $\pm$ 3.1\% & 12.2  & 91    & 0.2 $\pm$ 0.4 & 2.0\% $\pm$ 4.4\% & 0.1   & 82 \\
    \textbf{} & \textbf{} & \textbf{20\%} & 0.2 $\pm$ 0.4 & 2.4\% $\pm$ 4.7\% & 0.0     & 78    & 0.1 $\pm$ 0.3 & 1.0\% $\pm$ 3.4\% & 10.7  & 91    & 0.2 $\pm$ 0.4 & 2.4\% $\pm$ 4.7\% & 0.1   & 78 \\
    \textbf{} & \multirow{2}[2]{*}{\textbf{U[t,5t]}} & \textbf{10\%} & 0.3 $\pm$ 0.5 & 3.2\% $\pm$ 5.1\% & 0.0     & 70    & 0.2 $\pm$ 0.4 & 1.7\% $\pm$ 4.1\% & 13.4  & 85    & 0.3 $\pm$ 0.5 & 3.0\% $\pm$ 5.0\% & 0.1   & 72 \\
    \textbf{} & \textbf{} & \textbf{20\%} & 0.3 $\pm$ 0.5 & 3.4\% $\pm$ 5.2\% & 0.0     & 69    & 0.2 $\pm$ 0.4 & 2.1\% $\pm$ 4.4\% & 12.1  & 81    & 0.3 $\pm$ 0.5 & 3.4\% $\pm$ 5.2\% & 0.1   & 69 \\ \hline
    \multirow{3}[4]{*}{\textbf{3}} & \multirow{2}[2]{*}{\textbf{U[t,2t]}} & \textbf{10\%} & 0.3 $\pm$ 0.5 & 3.3\% $\pm$ 5.2\% & 0.0     & 69    & 0.1 $\pm$ 0.3 & 1.3\% $\pm$ 3.4\% & 16.4  & 86    & 0.3 $\pm$ 0.5 & 3.1\% $\pm$ 5.0\% & 0.1   & 71 \\
    \textbf{} & \textbf{} & \textbf{20\%} & 0.3 $\pm$ 0.4 & 2.8\% $\pm$ 5.0\% & 0.0     & 74    & 0.1 $\pm$ 0.3 & 0.9\% $\pm$ 3.0\% & 14.1  & 91    & 0.2 $\pm$ 0.4 & 2.4\% $\pm$ 4.7\% & 0.1   & 78 \\
    \textbf{} & \multirow{2}[2]{*}{\textbf{U[t,5t]}} & \textbf{10\%} & 0.4 $\pm$ 0.5 & 3.9\% $\pm$ 5.2\% & 0.0     & 62    & 0.2 $\pm$ 0.4 & 1.4\% $\pm$ 3.5\% & 17.0    & 85    & 0.4 $\pm$ 0.5 & 3.9\% $\pm$ 5.2\% & 0.0     & 62 \\
    \textbf{} & \textbf{} & \textbf{20\%} & 0.4 $\pm$ 0.5 & 3.9\% $\pm$ 5.1\% & 0.0     & 61    & 0.2 $\pm$ 0.4 & 1.6\% $\pm$ 3.5\% & 14.4  & 83    & 0.4 $\pm$ 0.5 & 3.9\% $\pm$ 5.1\% & 0.0     & 61 \\ \hline
    \multirow{3}[4]{*}{\textbf{4}} & \multirow{2}[2]{*}{\textbf{U[t,2t]}} & \textbf{10\%} & 0.4 $\pm$ 0.5 & 3.9\% $\pm$ 5.4\% & 0.0     & 63    & 0.2 $\pm$ 0.4 & 2.0\% $\pm$ 4.5\% & 16.6  & 82    & 0.4 $\pm$ 0.5 & 3.7\% $\pm$ 5.3\% & 0.1   & 65 \\
    \textbf{} & \textbf{} & \textbf{20\%} & 0.3 $\pm$ 0.5 & 3.5\% $\pm$ 5.1\% & 0.0     & 66    & 0.1 $\pm$ 0.3 & 1.3\% $\pm$ 3.7\% & 13.4  & 84    & 0.3 $\pm$ 0.5 & 3.4\% $\pm$ 5.1\% & 0.0     & 67 \\
    \textbf{} & \multirow{2}[2]{*}{\textbf{U[t,5t]}} & \textbf{10\%} & 0.5 $\pm$ 0.5 & 5.2\% $\pm$ 5.1\% & 0.0     & 46    & 0.3 $\pm$ 0.4 & 2.3\% $\pm$ 4.2\% & 18.7  & 75    & 0.5 $\pm$ 0.5 & 5.1\% $\pm$ 5.1\% & 0.0     & 47 \\
    \textbf{} & \textbf{} & \textbf{20\%} & 0.6 $\pm$ 0.5 & 5.2\% $\pm$ 5.1\% & 0.0     & 46    & 0.3 $\pm$ 0.5 & 2.6\% $\pm$ 4.3\% & 14.6  & 72    & 0.5 $\pm$ 0.5 & 5.2\% $\pm$ 5.0\% & 0.0     & 46 \\
    \bottomrule
    \end{tabular}%
  \label{tab:heuristic_cih_small_inst}%
\end{sidewaystable}%
%%%%%%%%%%%%%%%%%%%%%%%%%%%%%%%%%%%%%%%%%%%%%%%%%%%%%%%%%%%%%%%

The results of Table \ref{tab:heuristic_cih_small_inst} show the robust behavior of the proposed heuristic. Indeed, even in the most difficult situations (U[t,5t], Inc 20\%) the increase of 
stations is 5,2\% in average and the number of ties is 46. The computational times are practically null, which 
shows the computational efficiency of the heuristic. Concerning the CIH approaches with post-optimization
routines, we see that the greater flexibility of the LS1 causes significant improves in the
results of CIH, increasing the number of ties compared with the solutions of ALWIBP-1 model
in 24\% (272 instances). We also outline that, in some cases, the average number of additional stations 
in the assembly lines is reduced by half. Despite the fact that LS2 improves the quality of solutions in only a few instances,
it is very computationally efficient, running in very little time (0.1s in average).

Table \ref{tab:heuristic_cih_large_inst} shows the results obtained for the large instances. 
We notice that the number of Ties (66\%) is practically the same 
as observed in the previous case, but the average number of additional workstations is smaller. This states that the good results
of this scenario is powered by the flexibility of assembly lines with many tasks.

The computational times of the three approaches of the CIH show that the method is scalable. Moreover,
the LS1 increases the number of Ties related to the solutions of the ALWIBP-1 model, spending
up to 53s on average, a very encouraging result if we consider that the CPLEX exceeds the
computational time limit in most of cases.

% Table generated by Excel2LaTeX from sheet 'Heuristicas_ALWIBP_Resumo' %%%%%%%%%%%%%%%%%%%%%%%%%%
\begin{sidewaystable}[htbp]
  \scriptsize
  \centering
    \caption{Computational results: CIH approaches ($|N|$ = 100)}
    \begin{tabular}{ccc|cccc|cccc|cccc}
    \toprule
    \multirow{2}[0]{*}{\textbf{|W|}} & \multirow{2}[0]{*}{\textbf{Var}} & \multirow{2}[0]{*}{\textbf{Inc}} & \multicolumn{4}{c}{\textbf{CIH}} & \multicolumn{4}{c}{\textbf{CIH-LS1}} & \multicolumn{4}{c}{\textbf{CIH-LS2}} \\
    \cline{4-15}
          &       &       & $\mathbf{m_{h_{\uparrow}}} \pm \mathbf{\sigma_{m_{h_{\uparrow}}}}$ & $\mathbf{m_{h_{\uparrow}}} \pm \mathbf{\sigma_{m_{h_{\uparrow}}}} (\%)$ & $\mathbf{t_{h}(s)}$ & \textbf{Ties} & $\mathbf{m_{h_{\uparrow}}} \pm \mathbf{\sigma_{m_{h_{\uparrow}}}}$ & $\mathbf{m_{h_{\uparrow}}} \pm \mathbf{\sigma_{m_{h_{\uparrow}}}} (\%)$ & $\mathbf{t_{h}(s)}$ & \textbf{Ties} & $\mathbf{m_{h_{\uparrow}}} \pm \mathbf{\sigma_{m_{h_{\uparrow}}}}$ & $\mathbf{m_{h_{\uparrow}}} \pm \mathbf{\sigma_{m_{h_{\uparrow}}}} (\%)$ & $\mathbf{t_{h}(s)}$ & \textbf{Ties} \\ \hline
     \multirow{3}[4]{*}{\textbf{1}} & \multirow{2}[2]{*}{\textbf{U[t,2t]}} & \textbf{10\%} & 0.3 $\pm$ 0.4 & 1.5\% $\pm$ 2.6\% & 0.0     & 73    & 0.1 $\pm$ 0.3 & 0.5\% $\pm$ 1.6\% & 47.5  & 89    & 0.2 $\pm$ 0.4 & 1.1\% $\pm$ 2.3\% & 5.2   & 79 \\
    \textbf{} & \textbf{} & \textbf{20\%} & 0.3 $\pm$ 0.5 & 1.8\% $\pm$ 2.9\% & 0.0     & 71    & 0.2 $\pm$ 0.4 & 0.9\% $\pm$ 2.3\% & 46.9  & 82    & 0.2 $\pm$ 0.4 & 1.3\% $\pm$ 2.6\% & 6.6   & 77 \\
    \textbf{} & \multirow{2}[2]{*}{\textbf{U[t,5t]}} & \textbf{10\%} & 0.4 $\pm$ 0.5 & 1.9\% $\pm$ 2.7\% & 0.0     & 64    & 0.2 $\pm$ 0.4 & 0.8\% $\pm$ 2.0\% & 43.4  & 81    & 0.4 $\pm$ 0.5 & 1.9\% $\pm$ 2.7\% & 2.7   & 65 \\
    \textbf{} & \textbf{} & \textbf{20\%} & 0.3 $\pm$ 0.4 & 1.4\% $\pm$ 2.5\% & 0.0     & 73    & 0.1 $\pm$ 0.3 & 0.4\% $\pm$ 1.7\% & 40.2  & 87    & 0.3 $\pm$ 0.4 & 1.3\% $\pm$ 2.4\% & 2.8   & 75 \\ \hline
    \multirow{3}[4]{*}{\textbf{2}} & \multirow{2}[2]{*}{\textbf{U[t,2t]}} & \textbf{10\%} & 0.3 $\pm$ 0.4 & 1.3\% $\pm$ 2.4\% & 0.0     & 70    & 0.2 $\pm$ 0.4 & 0.8\% $\pm$ 1.9\% & 48.1  & 80    & 0.2 $\pm$ 0.4 & 1.0\% $\pm$ 2.2\% & 2.3   & 76 \\
    \textbf{} & \textbf{} & \textbf{20\%} & 0.3 $\pm$ 0.4 & 1.3\% $\pm$ 2.4\% & 0.0     & 75    & 0.2 $\pm$ 0.4 & 0.8\% $\pm$ 2.0\% & 46.6  & 82    & 0.2 $\pm$ 0.4 & 1.2\% $\pm$ 2.3\% & 2.1   & 77 \\
    \textbf{} &  \multirow{2}[2]{*}{\textbf{U[t,5t]}} & \textbf{10\%} & 0.3 $\pm$ 0.5 & 1.7\% $\pm$ 2.7\% & 0.1   & 68    & 0.2 $\pm$ 0.4 & 0.9\% $\pm$ 2.2\% & 46.3  & 74    & 0.3 $\pm$ 0.5 & 1.6\% $\pm$ 2.6\% & 3.0     & 70 \\
    \textbf{} & \textbf{} & \textbf{20\%} & 0.3 $\pm$ 0.5 & 1.6\% $\pm$ 2.6\% & 0.0     & 72    & 0.2 $\pm$ 0.4 & 0.9\% $\pm$ 2.2\% & 43.9  & 78    & 0.3 $\pm$ 0.5 & 1.6\% $\pm$ 2.6\% & 1.3   & 72 \\ \hline
    \multirow{3}[4]{*}{\textbf{3}} & \multirow{2}[2]{*}{\textbf{U[t,2t]}} & \textbf{10\%} & 0.3 $\pm$ 0.5 & 1.5\% $\pm$ 2.6\% & 0.1   & 70    & 0.2 $\pm$ 0.4 & 0.9\% $\pm$ 2.1\% & 49.6  & 78    & 0.3 $\pm$ 0.4 & 1.4\% $\pm$ 2.5\% & 3.1   & 72 \\
    \textbf{} & \textbf{} & \textbf{20\%} & 0.3 $\pm$ 0.5 & 1.6\% $\pm$ 2.6\% & 0.1   & 71    & 0.2 $\pm$ 0.4 & 0.6\% $\pm$ 1.9\% & 46.9  & 82    & 0.3 $\pm$ 0.4 & 1.4\% $\pm$ 2.4\% & 2.5   & 74 \\
    \textbf{} & \multirow{2}[2]{*}{\textbf{U[t,5t]}} & \textbf{10\%} & 0.4 $\pm$ 0.5 & 2.1\% $\pm$ 2.9\% & 0.1   & 63    & 0.2 $\pm$ 0.4 & 0.9\% $\pm$ 2.1\% & 53.4  & 79    & 0.4 $\pm$ 0.5 & 2.1\% $\pm$ 2.9\% & 0.7   & 63 \\
    \textbf{} & \textbf{} & \textbf{20\%} & 0.4 $\pm$ 0.5 & 2.2\% $\pm$ 2.8\% & 0.1   & 57    & 0.2 $\pm$ 0.4 & 0.9\% $\pm$ 2.2\% & 46.9  & 75    & 0.4 $\pm$ 0.5 & 2.2\% $\pm$ 2.8\% & 1.2   & 58 \\ \hline
   \multirow{3}[4]{*}{\textbf{4}} & \multirow{2}[2]{*}{\textbf{U[t,2t]}} & \textbf{10\%} & 0.3 $\pm$ 0.5 & 1.6\% $\pm$ 2.5\% & 0.1   & 68    & 0.3 $\pm$ 0.4 & 1.0\% $\pm$ 2.4\% & 48.3  & 68    & 0.3 $\pm$ 0.5 & 1.4\% $\pm$ 2.4\% & 1.7   & 70 \\
    \textbf{} & \textbf{} & \textbf{20\%} & 0.4 $\pm$ 0.5 & 1.8\% $\pm$ 2.8\% & 0.1   & 63    & 0.2 $\pm$ 0.4 & 0.9\% $\pm$ 2.2\% & 48.9  & 74    & 0.3 $\pm$ 0.5 & 1.7\% $\pm$ 2.7\% & 2.0     & 66 \\
    \textbf{} & \multirow{2}[2]{*}{\textbf{U[t,5t]}} & \textbf{10\%} & 0.4 $\pm$ 0.5 & 2.3\% $\pm$ 2.8\% & 0.1   & 56    & 0.2 $\pm$ 0.4 & 0.7\% $\pm$ 2.2\% & 50.5  & 72    & 0.4 $\pm$ 0.5 & 2.2\% $\pm$ 2.8\% & 2.3   & 57 \\
    \textbf{} & \textbf{} & \textbf{20\%} & 0.5 $\pm$ 0.5 & 2.7\% $\pm$ 2.9\% & 0.1   & 52    & 0.3 $\pm$ 0.4 & 1.3\% $\pm$ 2.5\% & 49.6  & 70    & 0.5 $\pm$ 0.5 & 2.5\% $\pm$ 2.9\% & 0.9   & 55 \\
    \bottomrule
    \end{tabular}%
  \label{tab:heuristic_cih_large_inst}%
\end{sidewaystable}%
%%%%%%%%%%%%%%%%%%%%%%%%%%%%%%%%%%%%%%%%%%%%%%%%%%%%%%%%%%%%%%% 

\subsubsection{Experiment 3: Validation of heuristics}

In this section, we present results of the CIH taking into account the
largest instances of the ALWIBP benchmark. Also, we validate our
method by comparing it with a more straightforward substitution heuristic.

CPLEX fails in obtaining solutions for the  ``very large'' instances described earlier. Although these instances are more theoretical than
practical, computational experiments with this scenario are important in order to put the scalability of the solving methods to proof.
Thus, we highlight the ability of the CIH to solve problems of these dimensions in a 
reasonable computational times (3 min on average), as
shown in Table \ref{tab:heuristic_cih_very_large_inst}. Additionally it is remarkable that in 21\% of the instances,
the number of stations of a conventional assembly line did not change with the insertion
of disabled workers.

% Table generated by Excel2LaTeX from sheet 'Heuristicas_ALWIBP_Resumo'
\begin{table}[htbp]
\centering
\caption{Computational results: CIH ($|N|$ = 1000)}
\begin{tabular}{ccc|cccc}
\toprule
\multirow{2}[0]{*}{\textbf{|W|}} & \multirow{2}[0]{*}{\textbf{Var}} & \multirow{2}[0]{*}{\textbf{Inc}} & \multicolumn{4}{c}{\textbf{CIH}}  \\
\cline{4-7}
& & & $\mathbf{m_{h_{\uparrow}}} \pm \mathbf{\sigma_{m_{h_{\uparrow}}}}$ & $\mathbf{m_{h_{\uparrow}}} \pm \mathbf{\sigma_{m_{h_{\uparrow}}}} (\%)$ & $\mathbf{t_{h}(s)}$ & \textbf{Ties}$_{sp1}$ \\
\hline
\textbf{1} & \textbf{U[t,2t]} & \textbf{10\%} & 0.8 $\pm$ 0.9 & 0.5\% $\pm$ 0.7\% & 50.1 & 34 \\
\textbf{} & \textbf{} & \textbf{20\%} & 0.9 $\pm$ 0.9 & 0.5\% $\pm$ 0.7\% & 50.8 & 31 \\
\textbf{} & \textbf{U[t,5t]} & \textbf{10\%} & 0.9 $\pm$ 0.9 & 0.6\% $\pm$ 0.7\% & 50.4 & 31 \\
\textbf{} & \textbf{} & \textbf{20\%} & 1.0 $\pm$ 0.9 & 0.6\% $\pm$ 0.7\% & 50.2 & 29 \\ \hline
\textbf{2} & \textbf{U[t,2t]} & \textbf{10\%} & 1.1 $\pm$ 0.9 & 0.7\% $\pm$ 0.7\% & 106.6 & 25 \\
\textbf{} & \textbf{} & \textbf{20\%} & 1.1 $\pm$ 0.9 & 0.7\% $\pm$ 0.7\% & 105.9 & 24 \\
\textbf{} & \textbf{U[t,5t]} & \textbf{10\%} & 1.5 $\pm$ 1.0 & 1.0\% $\pm$ 0.8\% & 105.9 & 16 \\
\textbf{} & \textbf{} & \textbf{20\%} & 1.5 $\pm$ 1.1 & 1.0\% $\pm$ 0.8\% & 105.9 & 17 \\ \hline
\textbf{3} & \textbf{U[t,2t]} & \textbf{10\%} & 1.2 $\pm$ 1.0 & 0.8\% $\pm$ 0.8\% & 167.1 & 21 \\
\textbf{} & \textbf{} & \textbf{20\%} & 1.3 $\pm$ 1.0 & 0.8\% $\pm$ 0.7\% & 167.3 & 22 \\
\textbf{} & \textbf{U[t,5t]} & \textbf{10\%} & 1.8 $\pm$ 1.1 & 1.1\% $\pm$ 0.8\% & 164.4 & 14 \\
\textbf{} & \textbf{} & \textbf{20\%} & 1.9 $\pm$ 1.0 & 1.2\% $\pm$ 0.8\% & 163.6 & 8 \\ \hline
\textbf{4} & \textbf{U[t,2t]} & \textbf{10\%} & 1.4 $\pm$ 1.0 & 0.9\% $\pm$ 0.8\% & 223.4 & 21 \\
\textbf{} & \textbf{} & \textbf{20\%} & 1.4 $\pm$ 1.0 & 0.9\% $\pm$ 0.8\% & 225.4 & 21 \\
\textbf{} & \textbf{U[t,5t]} & \textbf{10\%} & 2.1 $\pm$ 1.2 & 1.4\% $\pm$ 0.9\% & 224.1 & 10 \\
\textbf{} & \textbf{} & \textbf{20\%} & 2.1 $\pm$ 1.2 & 1.4\% $\pm$ 0.9\% & 222.0 & 9 \\
\bottomrule
\end{tabular}%
\label{tab:heuristic_cih_very_large_inst}%
\end{table}%

\newpage
In order to study the effectiveness of our algorithm in comparison to more straightforward strategies,  we have 
implemented a simple substitution heuristic (denoted by SH) and evaluated the obtained results. The idea behind the SH is simply to substitute conventional workers for disabled ones, maintaining the task assignments. Our goal with this experiment is to show that simple strategies such as this one, which lack the ability to perform more structural changes in the line structure are ineffective. The SH works as follows: given a SALBP-1
solution the algorithm searches for the most suitable stations in which to insert the disabled workers without 
changing task assignment and respecting cycle time and
task/worker feasibilities constraints. Algorithm \ref{alg:sh} presents the pseudo-code of
this heuristic. Note that the \textit{LD} formulation, modelled by equations
(\ref{eq:ld1})-(\ref{eq:ld7}), is part of the kernel of the SH.

\begin{algorithm} \caption{Simple Heuristic}
  \small
  \label{alg:sh}
  \begin{algorithmic}[1]
  \REQUIRE $x^{ref}$ (SALBP-1 reference solution);
  \STATE $\overline{m} \leftarrow$ \textit{stations}($x^{ref}$);  
  \STATE Let $M$ be an integer matrix with $\overline{m}$ rows and $|W|$ columns;
  \FOR{$s=1,...,\overline{m}$}
    \FORALL{$w \in W$}    
      \STATE $M_{sw} \leftarrow$ \textit{load}($s,w,x^{ref}$);
    \ENDFOR
  \ENDFOR
  \STATE $(l^*, y^*) \leftarrow$ \textit{solveLD}($M$);
  \IF{$l^* = \overline{C}$}
    \STATE \textbf{return} $y^*$;
  \ELSE
    \STATE \textit{there is no feasible solution};
  \ENDIF
  \end{algorithmic}
\end{algorithm}

\newpage
\begin{equation}
\label{eq:ld1}
\textrm{Min }   l
\end{equation}
\hspace{1cm} subject to
\begin{alignat}{2}
\label{eq:ld2}
l \geq \overline{C}, \\
\label{eq:ld3}
l \geq M_{sw}y_{sw}, && \qquad &s=1,...,\overline{m}, \forall w \in W, \\
\label{eq:ld4}
\sum_{s=1}^{\overline{m}} y_{sw} = 1, && \qquad &\forall w \in W, \\
\label{eq:ld5}
\sum_{w \in W} y_{sw} \leq 1, && \qquad &s=1,...,\overline{m}, \\
\label{eq:ld6}
y_{sw} \in \{0,1\}, && \qquad &s=1,...,\overline{m}, \forall w \in W, \\
\label{eq:ld7}
l \in \mathbb{Z}.
\end{alignat}

In lines 1-6, we compute the elements of the matrix $M$ as the load of
station $s$ if the ``conventional'' worker in that station is replaced 
by disabled worker $w$. To do that, we consider the \textit{load} function
(line 5). Next, in line 7, we solve the \textit{LD} formulation.
Lines 9-13 test the effectiveness of the SH. Thus, if the resulted
solution is feasible, that is, the best value found by \textit{LD}
model is equal to the desired cycle time, the algorithm returns
the variables $y$ which indicate the placement of each worker along
the assembly line. 

Computational experiments were 
conducted using the same input solutions as those used for the CIH, as well as the same machine settings
and CPLEX version. The SH could only find feasible solutions for about 4\% of the instances (3\%, 2\% and 7\% for the medium, large and very large instances, respectively), confirming our hypothesis that a more elaborate strategy such as the proposed CIH is necessary when dealing with the problem.

\section{Conclusions}

We propose the Assembly Line Worker Integration and Balancing Problem (ALWIBP), a new assembly line balancing problem arising in lines with conventional and disabled workers. This problem is relevant in a context where companies are urged to integrate disabled workers in their conventional productive schemes in order to cope with legislation issues or to include corporate social responsibility goals in the production planning process. 

To solve this problem, we first develop an integer linear model that minimizes the number of stations while ensuring
the presence of all disabled workers in the assembly line. From this model, one variant 
that reduces the idle time of stations with disabled workers was 
proposed.
We implement a heuristic approach called CIH, that starts with a simple
assembly line balancing situation and inserts the available disabled workers while reducing
additional stations. Furthermore, aiming the improvement of CIH solutions, two post-optimizations 
procedures based on MIP neighborhoods were designed. Finally, to validate heuristic results,
we also implemented a simple substitution heuristic (SH) that tries to insert disabled workers
in a conventional assembly line without changing its original task 
settings.

Results of an experimental study on a extensive number of instances indicate
the efficiency of the proposed heuristics, leading to the conclusion that
not only disabled workers can be included in the assembly lines with little productivity 
loss, but also that other planning goals can be simultaneously considered. Further work on 
this topic include the proposal of new adjacent objectives, implementation
of more sophisticated methods, and extensions
that cope with job rotation schemes, U-shaped assembly lines and some
lexicographic objectives.

\section{Acknowledgements}

This research was supported by CAPES-Brazil and MEC-Spain (coordinated project CAPES DGU 
258-12/PHB2011-0012-PC) and by FAPESP-Brazil. We also thank the project ``CORSARI
MAGIC DPI 2010-18243'' of the Ministerio de Ciencia e Innovación del Gobierno de España
within the Program ``Proyectos de Investigación Fundamental No Orientada''.
The authors also thank two anonymous reviewers and the editor for their 
useful comments.

%\bibliography{../../../bib/Claio2012}
\bibliography{library}

 \end{document}